\documentclass[preprint,12pt,compress]{elsarticle}

\usepackage{amssymb}
\usepackage{amsmath}
\usepackage[a4paper, body={16.8cm,23cm}]{geometry}

\begin{document}
\journal{arXiv}
\parindent 15pt
\parskip 8pt

\renewcommand{\arraystretch}{1.5}

 \newcommand{\eps}{\varepsilon}
 \newcommand{\lam}{\lambda}
 \newcommand{\To}{\rightarrow}
 \newcommand{\as}{{\rm d}P\times {\rm d}t-a.e.}
 \newcommand{\ps}{{\rm d}P-a.s.}
 \newcommand{\jf}{\int_t^T}
 \newcommand{\tim}{\times}

 \newcommand{\F}{\mathcal{F}}
 \newcommand{\E}{\mathbb{E}}
 \newcommand{\N}{\mathbb{N}}
 \newcommand{\s}{\mathcal{S}}
 \newcommand{\M}{{\rm M}}
 \newcommand{\T}{[0,T]}
 \newcommand{\LT}{L^2(\mathcal{F}_T;\R^k)}

 \newcommand{\R}{{\mathbb R}}
 \newcommand{\Q}{{\mathbb Q}}
 \newcommand{\RE}{\forall}

\newcommand{\gET}[1]{\underline{{\mathcal {E}}_{t,T}^g}[#1]}
\newcommand {\Lim}{\lim\limits_{n\rightarrow \infty}}
\newcommand {\Limk}{\lim\limits_{k\rightarrow \infty}}
\newcommand {\Limm}{\lim\limits_{m\rightarrow \infty}}
\newcommand {\llim}{\lim\limits_{\overline{n\rightarrow \infty}}}
\newcommand {\slim}{\overline{\lim\limits_{n\rightarrow \infty}}}
\newcommand {\Dis}{\displaystyle}

\begin{frontmatter}
\title {{\bfseries A general existence and uniqueness result on multidimensional BSDEs \tnoteref{found}}}
\tnotetext[fund]{Supported by the National Natural Science Foundation of China (No. 11371362), the China Postdoctoral Science Foundation (No. 2013M530173), the Qing Lan Project and the Fundamental Research Funds for the Central Universities (No. 2013RC20).\vspace{0.2cm}}

\author{ShaoYa XU\corref{cor}$^1$}

\author{ShengJun FAN\corref{cor1}$^{1,2}$}

\cortext[cor1]{Corresponding author. {\it Email Addresses:} f$\_$s$\_$j@126.com.\vspace{-0.2cm}}

\address{1. College of Sciences, China University of Mining and
Technology, Xuzhou 221116, PR China\vspace{0.1cm}\\
2. School of Mathematical Sciences, Fudan University, Shanghai 200433, PR China\vspace{-0.8cm}}

\begin{abstract}
This paper establishes a new existence and uniqueness result of
solutions for multidimensional backward stochastic differential
equations (BSDEs) whose generators satisfy a weak monotonicity condition and a general growth condition in $y$, which generalizes the
corresponding results in [2], [3] and
[5].\vspace{0.1cm}
\end{abstract}

\begin{keyword}
Backward stochastic differential equation \sep Existence and
uniqueness\sep Weakly monotonic condition\sep Lipschitz
condition\sep Mao's condition \vspace{0.2cm}

\MSC[2010] 60H10
\end{keyword}
\end{frontmatter}\vspace{-0.4cm}


\section{Introduction}

In this paper, we are concerned with the following multidimensional backward
stochastic differential equation (BSDE for short in the
remaining):
\begin{equation}
   y_t=\xi+\int_t^Tg(s,y_s,z_s){\rm d}s-\int_t^Tz_s {\rm d}B_s,\ \
    t\in\T,
\end{equation}
where $T>0$ is a constant called the time horizon, $\xi$ is a
$k$-dimensional random vector called the terminal condition, the
random function $g(\omega,t,y,z):\Omega\tim \T\tim {\R}^{k }\tim
{\R}^{k\times d}\To {\R}^k$ is progressively measurable for each
$(y,z)$, called the generator of BSDE (1), and $B$ is a
$d$-dimensional Brownian motion. BSDE (1) is denoted by BSDE
($\xi,T,g$). The solution ($y_{\cdot},z_{\cdot}$) is a pair of
adapted processes.

BSDEs were initially introduced in a nonlinear form in 1990 by
Pardoux and Peng [4], who established an existence and uniqueness result
for the adapted and squared integrable solutions of BSDEs under the Lipschitz assumption of the
generator $g$. From then on, many researchers have been working on this subject, and many applications have been found in mathematical finance, stochastic control, and partial differential equations, etc. In particular, an interesting and important question is how to improve the existence and uniqueness result of [4] by weakening the Lipschitz continuity condition on the generator $g$. Here, we would like to cite some efforts devoted to this direction and related closely to this paper. In 1995, Mao [3] obtained an existence and uniqueness result of a
solution for (1) where $g$ satisfies some kind of non-Lipschitz
condition in $y$ called usually the Mao's condition. In 1999, Pardoux [5] established an existence and
uniqueness result of a solution for (1) where $g$ satisfies some
kind of monotonicity condition and a general growth condition in $y$. Furthermore, in 2003, using the same monotonicity condition as in [5] and a more general growth condition in $y$ for $g$, Briand et al. [1] investigated the existence and uniqueness of a
solution for (1). Recently, under the general growth condition employed in [5] as well as a weaker monotonicity condition in $y$ for $g$, Fan and Jiang [2] proved an existence and uniqueness result of a solution for (1), which unifies the results obtained in [3] and [5].

The objective of this paper is to further generalize the existence and uniqueness result obtained in [2]. We establish a new existence and uniqueness result for
solutions of multidimensional BSDEs whose generators satisfy the
weaker monotonicity condition in $y$ put forward by [2] and the more general growth condition in $y$ employed in [1] (see
Theorem 1 in Section 3), which generalizes the corresponding results
in [3], [5] and [2]. Particularly, it should be mentioned that the integrability condition on the process ${\{ g(t,0,0)\} _{t \in [0,T]}}$ used in [2] is also weakened in Theorem 1 of this paper. The remainder is organized as follows. We introduce some preliminaries
and establish a technical proposition in Section 2, and put forward and prove
our main result in Section 3.

\section{ Preliminaries}

Let us fix a
number $T> 0$, and two positive integers $k$ and $d$. Let
$(\Omega,\F,P)$ be a probability space carrying a standard
$d$-dimensional Brownian motion $(B_t)_{t\geq 0}$. Let
$(\F_t)_{t\geq 0}$ be the natural $\sigma$-algebra generated by
$(B_t)_{t\geq 0}$ and $\F=\F_T$. In this paper, the Euclidean norm
of a vector $y\in\R^k$ will be defined by $|y|$, and for an $k\times
d$ matrix $z$, we define $|z|=\sqrt{{\rm Tr}(zz^*)}$, where $z^*$ is
the transpose of $z$. Let $\langle x,y\rangle$ represent the inner
product of $x,y\in\R^k$. We denote by $\LT$ the set of all
$\R^k$-valued, square integral and $\F_T$-measurable random vectors.
Let ${\s}^2(0,T;\R^k)$ denote the set of $\R^k$-valued, adapted and
continuous processes $(\phi_t)_{t\in\T}$ such that\vspace{0.22cm}
$$\|\phi\|_{{\s}^2}^2:=\E[\sup_{t\in\T} |\phi_t|^2]<+\infty.
\vspace{0.05cm}$$ Moreover, let ${\rm M}^2(0,T;\R^{k\times d})$
denote the set of $(\F_t)$-progressively measurable ${\R}^{k\times
d}$-valued processes $(\varphi _t)_{ t\in\T}$ such
that
$$\quad\|\varphi\|_{{\rm M}^2}^2:=\E\left[\int_0^T |\varphi_t|^2\
{\rm d}t\right]<+\infty. \vspace{0.1cm}$$
Obviously,
${\s}^2(0,T;\R^{k})$ is a Banach space and ${\rm
M}^2(0,T;\R^{k\times d})$ is a Hilbert space.

 As mentioned in the introduction, we will deal only with BSDEs which
are equations of type (1), where the terminal condition $\xi\in
\LT$, and the generator $g$ is $(\F_t)$-progressively measurable for
each $(y,z)$.

{\bf Definition 1} \ A pair of processes $(y_t,z_t)_{t\in\T}$ is
called a solution to BSDE (1), if $(y_t,z_t)_{t\in\T}\in
{\s}^2(0,T;\R^k)\times {\rm M}^2(0,T;\R^{k\times d})$ and satisfies (1).\vspace{0.1cm}

Now, let us introduce the following Proposition 1, which will play
an important role in the proof of our main result. In stating it, the following assumption on the generator
$g$ is useful:\vspace{0.2cm}

{(A)}\ \ $\as,\RE\ (y,z)\in\R^k\times\R^{k\times d}$, \ $\langle y, g(\omega,t,y,z)\rangle\leq
\psi(|y|^2)+\lambda|y||z|+|y|f_t,$\vspace{0.4cm}\\
where $\lambda>0$ is a constant, $(f_t)_{t \in [0,T]}$ is a nonnegative and $(\F_t)$-measurable
process with \\
$$\E\left[\left(\int_0^T f_t{\rm d}t\right)^2\right]
<+\infty,$$\\
and $\psi(\cdot)$ is a nondecreasing and concave function
from $\R^+$ to itself with $\psi(0)=0$.\vspace{0.1cm}

{\bf Proposition 1} \  Let $g$ satisfy (A) and $(y_t,z_t)_{t\in\T}$
be a solution to BSDE $(\xi,T,g)$. Then there exists a constant
$C>0$ depending only on $\lambda$ and $T$ such that for each $0\leq
u\leq t\leq T$,
$$
\begin{array}{lll}
&&\Dis \E\left[\left.\sup\limits_{r\in
[t,T]}|y_r|^2\right|\F_u\right]+\E\left[\left.\int_t^T
|z_s|^2\ {\rm d}s\right|\F_u\right]\\
&\leq & \Dis
C\left\{\E\left[\left.|\xi|^2\right|\F_u\right]+\int_t^T
\psi\left(\E\left[\left.|y_s|^2\right|\F_u\right]\right)\ {\rm
d}s+\E\left[\left.\left(\int_t^T f_s\ {\rm
d}s\right)^2\right|\F_u\right]\right\}.
\end{array}
$$
{\bf Proof.}\ Applying It\^{o}'s formula to $|y_t|^2$ leads that for
each $t\in \T$,
\begin{equation}
|y_t|^2+\int_t^T |z_s|^2\ {\rm d}s=|\xi|^2+2\int_t^T \langle
y_s,g(s,y_s,z_s)\rangle \ {\rm d}s-2\int_t^T\langle y_s,z_s{\rm
d}B_s\rangle.
\end{equation}
By assumption (A) and the inequality $ 2ab\leq 2a^2+{b^2/2}$ we have
\begin{equation}
\begin{array}{lll}
2\langle y_s,g(s,y_s,z_s)\rangle &\leq &2\psi(|y_s|^2)+2
\lambda|y_s||z_s|+2|y_s|f_s\\
 & \leq &
 2\psi(|y_s|^2)+2\lambda^2|y_s|^2+{1\over 2}|z_s|^2+2|y_s|f_s.
\end{array}
\end{equation}
It follows from the Burkholder-Davis-Gundy inequality that $\{M_t:=\int_0^t\langle
y_s,z_s{\rm d}B_s\rangle\}_{t\in\T}$ is a uniformly integrable
martingale. In fact, for each $0\leq u\leq t\leq T$, we have
\begin{equation}
\begin{array}{lll}
\Dis 2\E\left[\left. \sup\limits_{r\in
[t,T]}\left|\int_r^T\langle y_s,z_s{\rm
d}B_s\rangle\right|\right|\F_u\right]&\leq & \Dis
2c\E\left[\left.\sup\limits_{r\in [t,T]}|y_r|\cdot\left(
\int_t^T|z_s|^2\ {\rm
d}s\right)^{1/2}\right|\F_u\right]\\
&\leq &\Dis {1\over 2}\E\left[\left.\sup\limits_{r\in
[t,T]}|y_r|^2\right|\F_u\right]+2c^2\E\left[\left.\int_t^T|z_s|^2\
{\rm d}s\right|\F_u\right]\\
&<&+\infty,
\end{array}
\end{equation}
where $c>0$ is a constant. Then, it follows from (2), (3) and (4) that for each $0\leq u\leq t\leq T$,
\begin{equation}
\Dis {1\over 2}\E\left[\left.\int_t^T |z_s|^2\ {\rm
d}s\right|\F_u\right]\leq \E\left[\left.X_t\right|\F_u\right]
+2\E\left[\left.\int_t^T|y_s|f_s{\rm d}s\right|\F_u\right],
\end{equation}
where
$$X_t=|\xi|^2+2\lambda^2\int_t^T |y_s|^2{\rm d}s+2\int_t^T
\psi(|y_s|^2)\  {\rm d}s.$$

Furthermore, by virtue of  (3), (4) and the following inequality
\begin{equation}
\begin{array}{lll}
\Dis 2\E\left[\left.\int_t^T|y_s|f_s{\rm d}s\right|\F_u\right]& \leq
&\Dis 2\E\left[\left.\sup_{r\in [t,T]}|y_r|\cdot\int_t^T f_s{\rm
d}s\right|\F_u\right]\\
&\leq & \Dis {1\over 4}\E\left[\left.\sup_{r\in
[t,T]}|y_r|^2\right|\F_u\right]+4\E\left[\left.\left(\int_t^T f_s\
{\rm d}s\right)^2\right|\F_u\right],
\end{array}
\end{equation}
it follows from (2) that for each $0\leq u\leq t\leq
T$,
$$
\begin{array}{lll}
&& \Dis {1\over 4}\E\left[\left.\sup\limits_{r\in
[t,T]}|y_r|^2\right|\F_u\right]+{1\over 2}\E\left[\left.\int_t^T |z_s|^2
\ {\rm d}s\right|\F_u\right]\\
&\leq & \Dis
\E\left[\left.X_t\right|\F_u\right]+4\E\left[\left.\left(\int_t^T
f_s\ {\rm d}s\right)^2\right|\F_u\right]
\Dis+2c^2\E\left[\left.\int_t^T |z_s|^2\ {\rm d}s\right|\F_u\right].
\end{array}
$$ Combining the above inequality, (5) and (6) with
$4$ being replaced by $32c^2$ yields that for each $0\leq u\leq
t\leq T$,
$$\begin{array}{lll}
&& \Dis \hspace{-1.8cm}\quad{1\over 8}\E\left[\left.\sup\limits_{r\in
[t,T]}|y_r|^2\right|\F_u\right]+{1\over 2}\E\left[\left.\int_t^T |z_s|^2\ {\rm
d}s\right|\F_u\right]\\
\hspace{-1.5cm}\leq&& \Dis \hspace{-1.6cm}(4c^2+1)\E\left[\left.X_t\right|\F_u\right]
+(16c^2+4)\E\left[\left.\left(\int_t^T f_s\ {\rm
d}s\right)^2\right|\F_u\right],
\end{array}$$ and then, in view of the definition
of $X_t$, Fubini's theorem, the concavity of $\psi(\cdot)$ and
Jensen's inequality, we have
$$\begin{array}{lll}
\quad\quad&& \Dis{1\over 8}\E\left[\left.\sup\limits_{r\in
[t,T]}|y_r|^2\right|\F_u\right]+{1\over 2}\E\left[\left.\int_t^T |z_s|^2\ {\rm
d}s\right|\F_u\right]\\
\quad\quad\quad\quad&\leq & \Dis
(4c^2+1)\E\left[\left.|\xi|^2\right|\F_u\right]+2(4c^2+1)\int_t^T
\psi\left(\E\left[\left.|y_s|^2\right|\F_u\right]\right)\ {\rm
d}s\\
&& \Dis +(16c^2+4)\E\left[\left.\left(\int_t^T f_s\ {\rm
d}s\right)^2\right|\F_u\right]+2\lambda^2(4c^2+1)\int_t^T
\E\left[\left.\sup\limits_{r\in [s,T]}|y_r|^2\right|\F_u\right]{\rm
d}s,
\end{array}$$
from which together with Gronwall's inequality, the desired result follows. The
proof is then completed.

{\bf Remark 1}\quad Proposition 1 improves the corresponding result in
[2], where the process $(f_t)$ defined in assumption (A)
is assumed to satisfy the condition that
$$\E\left[\int_0^T |f_t|^2{\rm
d}t\right]<+\infty.$$

\section{ Main result and its proof }

In this section, we will put forward and prove our main result. Let us first introduce the following
assumptions on the generator $g$:

(H1) \ $g$ satisfies the weakly monotonic condition in $y$, i.e.,
there exists a nondecreasing and concave function
$\kappa(\cdot):\R^+\mapsto \R^+$ with $\kappa(0)=0$, $\kappa(u)>0$
for $u>0$ and $\int_{0^+} {{\rm d}u\over \kappa(u)}=+\infty$ such
that $\as$,  $$\RE y_1,y_2\in \R^k,z\in\R^{k\times d}, \ \
\langle
y_1-y_2,g(\omega,t,y_1,z)-g(\omega,t,y_2,z)\rangle\leq
\kappa(|y_1-y_2|^2).$$

(H2) $\as$, $\RE\ z\in {\R^{k\times d}},\ \ \ y\longmapsto
g(\omega,t,y,z)$ is continuous.\vspace{0.3cm}

(H3) $\RE\ \alpha>0,\ \phi_\alpha(t):=\sup\limits_{|y|\leq \alpha}
|g(\omega,t,y,0)-g(\omega,t,0,0)|\in L^1(\T\tim \Omega)$.

(H4)\  $g$ is Lipschitz continuous in $z$ uniformly with respect to
$(\omega,t,y)$, i.e., there exists a
constant $\mu\geq 0$ such that $\as,$ \vspace{-0.1cm}
$$\RE\ y\in \R^k,z_1,z_2\in\R^{k\times d},\ \
|g(\omega,t,y,z_1)-g(\omega,t,y,z_2)|\leq \mu
|z_1-z_2|.$$

(H5)\  $\Dis\E\left[\left(\int_0^T |g(\omega,t,0,0)|\ {\rm
d}t\right)^2\right]<+\infty$.\\

In this paper, we want to obtain an existence and uniqueness result for BSDE
(1) under the previous assumptions (H1)-(H5) and $\xi\in\LT$. Firstly, let us recall a result in [2], which unifies the existence and uniqueness results obtained in [3] and [5]. For this, let us
introduce the following assumptions:

(H3') $g$ has a general growth with respect to $y$, i.e, $\as,$
$$\RE\ y\in \R^k,\ \ |g(\omega,t,y,0)|\leq
|g(\omega,t,0,0)|+\varphi(|y|),$$
where $\varphi:{\R}^+\To {\R}^+$ is an increasing continuous function.\vspace{0.2cm}

(H5')\  $\Dis\E\left[\int_0^T |g(\omega,t,0,0)|^2\ {\rm
d}t\right]<+\infty$.\vspace{0.2cm}

{\bf Proposition 2} (see Theorem 2.1 in \citep{Fan10})\ Let
assumptions (H1), (H2), (H3'), (H4) and (H5') hold. Then for each $\xi\in\LT$, BSDE $(\xi,T,g)$
has a unique solution.

The following Theorem 1 is the main result of this
paper.

{\bf Theorem 1}\ \  Let assumptions (H1)-(H5) hold. Then for each $\xi\in\LT$, BSDE $(\xi,T,g)$ has a unique
solution.

{\bf Remark 2}\ \ \ Note that (H3) and (H5) are strictly weaker than (H3') and (H5') respectively. It is clear that Theorem 1 generalizes Proposition 2 and the
corresponding results in [3, 5].

{\bf Example 1}\ \ Let $k = 2$ and for each $y=(y_1,y_2)\in\mathbb{R}^{2}$ and $z\in\mathbb{R}^{2\times d}$, let  $g(t,y,z) = ({g_1}(t,y,z),{g_2}(t,y,z))$ be defined by
$$
{g_i}(t,y,z) = |B_t|\cdot{e^{ - {y_i}}} + h(|y|) + |z| + \frac{1}{{\sqrt t }}\cdot {1_{t > 0}}, \ \ i=1,2,
$$
where
$$ h(x) = \left\{ {\begin{array}{*{20}{l}}
{ - x\ln x}\\
{{h^{'}}(x)(x - \delta ) + h(\delta )}\\
0
\end{array}} \right.\begin{array}{*{20}{c}}
,\\
,\\
,
\end{array}\begin{array}{*{20}{l}}
{0 < x \le \delta ;}\\
{x > \delta ;}\\
{\begin{array}{*{20}{c}}
{other}&{cases}
\end{array}}
\end{array}
$$
with $\delta>0$ small enough.

It is not hard to check that this $g$ satisfies (H3) and (H5), but does not satisfy (H3') and (H5'). At the same time, it is clear that $g$ satisfies (H2) and (H4) with $\mu=1$. In addition, we can also prove that $g$ satisfies (H1) (see Examples 2.4-2.5 and Remark 2.2 in [2] for details). Then, it follows from Theorem 1 that for each $\xi\in\LT$, BSDE$(\xi,T,g)$ has a unique solution. It should be mentioned that this conclusion can not be obtained by any known results including the previous proposition 2.\vspace{0.2cm}

{\bf The Proof of Theorem 1.} \ Assume that $g$ satisfies
assumptions (H1)-(H5). The proof of the uniqueness part is similar to that of the uniqueness part of Theorem 2.1 in [2], so we omit it.
Let us turn to the existence part. The proof will be split into two
steps.

First step: We shall prove that under assumptions (H1)-(H5),
provided that there exists a constant $K>0$ such that
\begin{equation}
\ps,\ |\xi|\leq K,\ \ {\rm and}\ \  \as,\ |g(t,0,0)|\leq K,
\end{equation}
BSDE $(\xi,T,g)$ has a solution.

For some large enough integer $\alpha>0$ which will be chosen later,
let $\theta_\alpha$ be a smooth function such that $0\leq
\theta_\alpha\leq 1$, $\theta_\alpha(y)=1$ for $|y|\leq \alpha$ and
$\theta_\alpha(y)=0$ as soon as $|y|\geq \alpha+1$. For each $n\geq
1$ and $z\in \R^{k\times d}$, we denote $q_n(z)=zn/(|z|\vee n)$ and
set\vspace{0.1cm}
$$h_n(t,y,z):=\theta_\alpha(y)(g(t,y,q_n(z))-g(t,0,0))
{n\over \phi_{\alpha+1}(t)\vee n}+g(t,0,0),$$
where $\phi_{\alpha}(\cdot)$ is defined in (H3).

It is clear from (7) that $h_n$ satisfies assumptions (H2) and (H4) and
(H5') for each $n\geq 1$. It is also easy from (7) and (H3) to check that
$|h_n(t,y,0)|\leq n+K$, which means that $h_n$ satisfies (H3'). We
now prove that $h_n$ satisfies also assumption (H1) but with
another concave function $\bar\kappa(\cdot)$ which will be chosen
later. Indeed, let us pick $y_1$ and $y_2$ in $\R^k$. If
$|y_1|>\alpha+1$ and $|y_2|>\alpha+1$, (H1) is trivially satisfied
and thus we reduce to the case where $|y_2|\leq \alpha+1$. We write
$$\begin{array}{lll}
&&\Dis \langle y_1-y_2,h_n(t,y_1,z)- h_n(t,y_2,z)\rangle\\
&=& \Dis \theta_\alpha(y_1){n\over \phi_{\alpha+1}(t)\vee n}\langle
y_1-y_2,g(t,y_1,q_n(z))- g(t,y_2,q_n(z))\rangle\\
&& +\Dis {n\over \phi_{\alpha+1}(t)\vee
n}(\theta_\alpha(y_1)-\theta_\alpha(y_2))\langle
y_1-y_2,g(t,y_2,q_n(z))- g(t,0,0)\rangle.
\end{array}
$$
Since $g$ satisfies (H1), the first term of the right-hand side of
the previous equality is smaller than the term
$\kappa(|y_1-y_2|^2)$. For the second term, we can use the fact
that $\theta_\alpha$ is $C(\alpha)$-Lipschitz, to get, since
$|y_2|\leq \alpha+1$,
$$\begin{array}{lll}
&&\Dis(\theta_\alpha(y_1)-\theta_\alpha(y_2))\langle
y_1-y_2,g(t,y_2,q_n(z))-
g(t,0,0)\rangle\\
&\leq & \Dis C(\alpha)|y_1-y_2|^2|g(t,y_2,q_n(z))-g(t,0,0)|\leq
C(\alpha)(\phi_{\alpha+1}(t)+\mu n)|y_1-y_2|^2
\end{array}$$
and thus
$${n\over \phi_{\alpha+1}(t)\vee
n}(\theta_\alpha(y_1)-\theta_\alpha(y_2))\langle
y_1-y_2,g(t,y_2,q_n(z))- g(t,0,0)\rangle\leq C(\alpha)(1+\mu
)n|y_1-y_2|^2.$$Hence, letting $\bar\kappa(x)=C(\alpha)(1+\mu
)nx+\kappa(x)$, we have
$$
\langle y_1-y_2,h_n(t,y_1,z)- h_n(t,y_2,z)\rangle\leq
\bar\kappa(|y_1-y_2|^2).
$$
It is clear that $\bar\kappa(\cdot)$ is a nondecreasing concave
function with $\bar\kappa(0)=0$ and $\bar\kappa(u)>0$ for $u>0$. Moreover, it follows from the concavity of $\kappa(\cdot)$ that
$$
\kappa(u)=\rho(u\cdot 1+(1-u)\cdot 0)\geq
u\kappa(1)+(1-u) \rho(0)=u\kappa(1),\ \ u\in [0,1],
$$ and then\vspace{0.1cm}
$$\int_{0^+} {{\rm d}u\over \bar\kappa(u)}=\int_{0^+} {{\rm
d}u\over C(\alpha)(1+\mu )nu+\kappa(u)}\geq {\kappa(1)\over
C(\alpha)(1+\mu )n+\kappa(1)}\int_{0^+} {{\rm d}u\over
\kappa(u)}=+\infty.\vspace{0.1cm}$$
Then the pair $(\xi,h_n)$
satisfies all the assumptions of Proposition 2. Hence, for each
$n\geq 1$, BSDE $(\xi,T,h_n)$ has a unique solution
$(y_t^n,z_t^n)_{t\in \T}$.

Furthermore, it follows from (H1), (H4) and (7) that
$$\begin{array}{lll}
\Dis \langle y, h_n(t,y,z)\rangle &=& \Dis \theta_\alpha(y){n\over
\phi_{\alpha+1}(t)\vee n}\langle y,
g(t,y,q_n(z))-g(t,0,q_n(z))\\
&& \Dis \quad\quad\quad\quad\quad\quad\quad\quad+g(t,0,q_n(z))-g(t,0,0)\rangle+\langle y,g(t,0,0)\rangle\\
&\leq & \kappa(|y|^2)+\mu |y||z|+K|y|.
\end{array}$$
Consequently, assumption (A) is satisfied for the generator
$h_n$ of BSDE $(\xi,T,h_n)$ with $\psi(u)=\kappa(u)$,
$\lambda=\mu$ and $f_t\equiv K$. It then follows from Proposition 1
and (7) that there exists a constant $C>0$ depending only on $\mu$
and $T$ such that for each $0\leq u\leq t\leq T$,
$$
\Dis \E\left[\left.\sup\limits_{r\in
[t,T]}|y_r^n|^2\right|\F_u\right]+\E\left[\left.\int_t^T
|z_s^n|^2\ {\rm d}s\right|\F_u\right]\\
\leq  \Dis CK^2(1+T^2)+C\int_t^T
\kappa\left(\E\left[\left.|y_s^n|^2\right|\F_u\right]\right)\ {\rm
d}s.
$$
Since $\kappa(\cdot)$ is a nondecreasing and concave function with
$\kappa(0)=0$, it increases at most linearly, i.e., there exists a
constant $A>0$ such that $\kappa(u)\leq A(u+1)$ for each $u\geq 0$.
Applying Gronwall's inequality to the previous inequality yields
that
$$
\E\left[\left.|y^n_t|^2\right|\F_u\right]+\E\left[\left.\int_t^T
|z^n_s|^2\ {\rm d}s\right|\F_u\right]\leq \alpha^2,$$ where
$\alpha:=\sqrt {CK^2(1+T^2)+CAT}\cdot e^{{CAT/2}}$. Substituting
$u=t$ in the previous inequality yields that for each $n\geq 1$ and
$t\in\T$,\vspace{0.1cm}
\begin{equation}
|y^n_t|\leq \alpha,\ \ {\rm and}\ \ \ \mathbb{E}\left[\int_0^T
|z^n_s|^2\ {\rm d}s\right]\leq \alpha^2.\vspace{0.1cm}
\end{equation}
As a byproduct, $(y^n_t,z^n_t)_{t\in \T}$ solves the BSDE
$(\xi,T,g_n)$, where\vspace{0.1cm}
$$g_n(t,y,z)=(g(t,y,q_n(z))-g(t,0,0)){n\over
\phi_{\alpha+1}(t)\vee n} +g(t,0,0).$$

In the sequel, for each $n\geq 1$ and $i\geq 1$, let
$\hat{y}_\cdot^{n,i}=y_\cdot^{n+i}-y_\cdot^{n}$,
$\hat{z}_\cdot^{n,i}=z_\cdot^{n+i}-z_\cdot^{n}$. We have
$$
\hat{y}_t^{n,i}=\int_t^T \hat{g}^{n,i}(s,\hat{y}_s^{n,i}, \hat
z_s^{n,i})\ {\rm d}s-\int_t^T \hat{z}_s^{n,i}{\rm d}B_s,\ \ \ t\in
\T, $$ where for each $y\in \R^k$,
$$\begin{array}{lll}
\hat{g}^{n,i}(s,y,z)&:= & \Dis
(g(s,y+y_s^n,q_{n+i}(z+z_s^n))-g(s,0,0)){(n+i)\over
\phi_{\alpha+1}(s)\vee (n+i)}\\
&& \Dis -(g(s,y_s^n,q_n(z_s^n))-g(s,0,0)){n\over
\phi_{\alpha+1}(s)\vee n}.
\end{array}$$
It also follows from (H1) and (H4) that
\begin{equation}
\begin {array}{lll}
\hspace{-0.6cm}$$\langle y,\hat{g}^{n,i}(s,y,z)\rangle &=& \Dis {(n+i)\over
\phi_{\alpha+1}(s)\vee
(n+i)} \langle y, g(s,y+y_s^n,q_{n+i}(z+z_s^n))-g(s,y_s^n,q_n(z_s^n))\rangle\\
&& \Dis +({(n+i)\over \phi_{\alpha+1}(s)\vee (n+i)}-{n\over
\phi_{\alpha+1}(s)\vee n})
\langle y,g(s,y_s^n,q_{n}(z_s^n))-g(s,0,0)\rangle \\
&\leq& \kappa(|y|^2)+\mu|y|(|z|+2|z_s^n|1_{|z_s^n|>n})
+21_{\phi_{\alpha+1}(s)>n}|y| (\phi_{\alpha+1}(s)+\mu |z_s^n|),$$
\end{array}
\end{equation}
where we have used the fact that
$$\begin{array}{lll}
|q_{n+i}(z+z_s^n)-q_n(z_s^n)| &\leq
&|q_{n+i}(z+z_s^n)-q_{n+i}(z_s^n)|
+|q_{n+i}(z_s^n)-q_n(z_s^n)|\\
&\leq & |z|+2|z_s^n|1_{|z_s^n|>n}.
\end{array}$$
Then, combining (8), (9) and the inequality $2ab\leq 2a^2+{b^2/2}$ we deduce that
$$\begin{array}{lll}
2\langle \hat{y}_s^{n,i},\hat{g}^{n,i}(s,\hat{y}_s^{n,i},\hat{z}_s^{n,i})\rangle &\leq &2\kappa(|\hat{y}_s^{n,i}|^2)+2
\mu^2|\hat{y}_s^{n,i}|^2+{1\over 2}|\hat{z}_s^{n,i}|^2+4\alpha\mu|z_s^n|1_{|z_s^n|>n}\\
&& \Dis +4\alpha1_{\phi_{\alpha+1}(s)>n}(\phi_{\alpha+1}(s)+\mu |z_s^n|).
\end{array}$$
With this inequality in hand, using a similar argument to the proof of Proposition 3.1 in [2], we can deduce that there
exists a constant $C>0$ depending only on $\mu$ and $T$ such that
for each $t\in \T$ and each $n,i\geq1$,
$$\begin{array}{lll}
&& \Dis \E\left[\sup\limits_{r\in [t,T]}|\hat{y}_r^{n,i}|^2 +\int_t^T
|\hat{z}_s^{n,i}|^2\ {\rm d}s\right]\\
&\leq & \Dis C\int_t^T \kappa\left(\E\left[\sup\limits_{r\in
[s,T]}|\hat{y}_r^{n,i}|^2\right]\right)\ {\rm
d}s+2C\alpha\mu\E\left[\int_t^T
|z_s^{n}|1_{|z_s^{n}|>n}\ {\rm d}s\right]\\
&& \Dis +2C\alpha\E\left[\int_t^T
1_{\phi_{\alpha+1}(s)>n}(\phi_{\alpha+1}(s)+\mu |z_s^{n}|)\ {\rm
d}s\right].
\end{array}$$
Furthermore, with the help of (8), (H3) and the assumptions of
$\kappa(\cdot)$, taking the limsup with respect to $n$ in the
previous inequality and using Fatou's lemma and Bihari's inequality
yields that $\{(y^n_t,z^n_t)_{t\in \T}\}_{n=1}^{\infty}$ is a Cauchy
sequence in the process space ${\s}^2(0,T;\R^k)\times {\rm
M}^2(0,T;\R^{k\times d})$. Finally, we can pass to the limit in the
approximating BSDE $(\xi,T,g_n)$, which yields a solution to BSDE
$(\xi,T,g)$.

Second step: We now treat the general case. For each $n\geq 1$, let
\begin{equation}
\xi_n:=q_n(\xi)\ \ {\rm and}\ \
g_n(t,y,z):=g(t,y,z)-g(t,0,0)+q_n(g(t,0,0)).
\end{equation}
Clearly, the $(\xi_n,g_n)$ satisfies the assumptions of the first
step and
\begin{equation}
\E\left[|\xi_n-\xi|^2\right]\To 0,\ \ \E\left[\left(\int_0^T
|q_n(g(s,0,0))-g(s,0,0)|\ {\rm d}s\right)^{2}\right]\To 0
\end{equation}
as $n\To \infty$ by (H5). For each $n\geq 1$, thanks to the first
step of this proof, let $(y^n_t,z^n_t)_{t\in \T}$ denote the unique
solution to BSDE $(\xi_n,T,g_n)$. For each $n\geq 1$ and $m\geq
1$, let $\hat{y}_\cdot^{n,m}=y_\cdot^{n}-y_\cdot^{m}$,
$\hat{z}_\cdot^{n,m}=z_\cdot^{n}-z_\cdot^{m}$, we have
\begin{equation}
\hat{y}_t^{n,m}=\xi_n-\xi_m+\int_t^T
\hat{g}^{n,m}(s,\hat{y}_s^{n,m}, \hat{z}_s^{n,m})\ {\rm d}s-\int_t^T
\hat{z}_s^{n,m}{\rm d}B_s,\ \ \ t\in \T,
\end{equation}
where for each $(y,z)\in \R^k\times \R^{k\times d}$,
$$\hat{g}^{n,m}(s,y,z):=g_{n}(s,y+y_s^m,z+z_s^m)-g_m(s,y_s^m,z_s^m).$$
We write
$$
\begin {array}{lll}
\langle y, \hat{g}^{n,m}(t,y,z)\rangle&=& \langle
y,g_{n}(t,y+y_t^m,z+z_t^m)
-g_{m}(t,y+y_t^m,z+z_t^m)\rangle\\
&& +\langle y, g_{m}(t,y+y_t^m,z+z_t^m)-g_m(t,y_t^m,z_t^m)\rangle.
\end{array}
$$
It follows from (10), (H1) and (H4) that for each $(y,z)\in
\R^k\times \R^{k\times d}$, $\as$,
$$
\begin {array}{lll}
\quad\langle y, \hat{g}^{n,m}(t,y,z)\rangle
&=&\langle y,q_n(g(t,0,0))-q_m(g(t,0,0))\rangle\\
&& +\langle y, g(t,y+y_t^m,z+z_t^m)-g(t,y_t^m,z_t^m)\rangle\\
&\leq & |y||q_n(g(t,0,0))-q_m(g(t,0,0))|+\kappa(|y|^2)+\mu
|y||z|.\end{array} $$
Consequently, assumption
(A) is satisfied for the generator $\hat g^{n,m}(t,y,z)$ of BSDE (12) with $\psi(u)=\kappa(u)$, $\lambda=\mu$ and $f_t=
|q_n(g(t,0,0))-q_m(g(t,0,0))|$.  It then follows from Proposition 1 with $u=0$
that there exists a constant $C>0$ depending only on
$T$ and $\mu$ such that for each $t\in [0,T]$,\vspace{0.3cm}
\begin{equation}
\begin{array}{lll}
\hspace{-7cm}&& \Dis\E\left[\sup\limits_{r\in [t,T]}|\hat
y^{n,m}_r|^2\right]+\E\left[\int_t^T |\hat z^{n,m}_s|^2\ {\rm
d}s\right]\\
&\leq & \Dis C\E\left[|\xi_n-\xi_m|^2\right]+C\int_t^T\kappa
\left(\E\left[\sup_{r\in
[s,T]}|\hat y^{n,m}_r|^2 \right]\right){\rm d}s\\
&& \Dis +C\E\left[\left(\int_0^T |q_n(g(s,0,0))-q_m(g(s,0,0))|\ {\rm
d}s\right)^2\right].
\end{array}
\end{equation}
Note that there exists a constant $A>0$ such that $\kappa(u)\leq
A(u+1)$ for each $u\geq 0$. Gronwall's inequality yields that
for each $t\in \T$ and each $n,m\geq 1$,
$$
\begin{array}{lll}
\quad\quad\quad\quad&&\Dis \E\left[\sup\limits_{r\in [t,T]}|\hat
y^{n,m}_r|^2\right]+\E\left[\int_t^T |\hat z^{n,m}_s|^2\ {\rm
d}s\right]\\
&\leq & \Dis e^{CAT}\cdot
\left(4C\E\left[|\xi|^2\right]+CAT+4C\E\left[\left(\int_0^T
|g(s,0,0)| {\rm d}s\right)^2\right]\right).
\end{array}
$$
Thus, in view of (11), by taking the limsup in (13) with respect to
$n,m$ and using Fatou's lemma, the monotonicity and continuity of
$\kappa(\cdot)$ and Bihari's inequality we know that
$\{(y^n_t,z^n_t)_{t\in \T}\}_{n=1}^{\infty}$ is a Cauchy sequence in
the process space ${\s}^2(0,T;\R^k)\times {\rm M}^2(0,T;\R^{k\times
d})$. Let $(y_t,z_t)_{t\in \T}$ be the limit process of the sequence
$\{(y^n_t,z^n_t)_{t\in \T}\}_{n=1}^{\infty}$. We pass to the limit
in uniform convergence in probability for BSDEs $(\xi_n,T,g_n)$, thanks to (H2), (H3) and (H4), to
see that $(y_t,z_t)_{t\in \T}$ solves BSDE $(\xi,T,g)$. Thus, we
complete the proof of Theorem 1.\vspace{0.4cm}




\end{document}